\documentclass[12pt]{article}
\usepackage{latexsym}
\usepackage{amsmath}
\usepackage{amssymb}
\usepackage{smacros2}
\usepackage{fullpage}
\usepackage{proof}
\usepackage{psfig}
\usepackage{meaning}
\usepackage{graphicx}
\usepackage[all]{xy}
\usepackage{mathabx}
\usepackage{xcolor}

\begin{document}
\def\ssp{\hspace*{0.3cm}}
\def\nlsp{\vspace{1ex}\\}
\def\mbf{\boldsymbol}
\def\und{\underline}
\def\limp{\mathrel{-\!\circ}}
\def\bull{\vrule height 1.3ex width 1.0ex depth -.1ex } %%black box
\def\ap{\raisebox{0.3ex}{\bf .}} 
\renewcommand{\phi}{\varphi}
\def\tsum{{\textstyle \sum}}
\def\tbig{{\textstyle \bigoplus}}
\def\tco{{\textstyle \coprod}}
\def\tbigu{{\textstyle \biguplus}}
\def\dquad{\quad\quad}
\def\mA{\mathbb A}
\def\mB{\mathbb B}
\def\mT{\mathbb T}
\def\mD{\mathbb D}
\def\mZ{\mathbb Z}
\def\mC{\mathbb C}
\def\vtube{\text{iff}}

\def\PAR{\mathrel{\invamp}} %%Use \mathrel to make a binary relation
\def\wp{\mathrel{\Invamp}} %%Looks the same but is faster to draw
\def\FR{\mathcal{\textbf{FR}}_k}
\def\UFR{\mathcal{\textbf{UFR}}_k}
\def\Ab{\textbf{Ab}}
\def\lca{\textbf{LCA}}
\def\mR{\mathbb{R}}
\def\im{\text{im}}

\title{A categorical approach to additive combinatorics}%\\Version 2.0}
%Journal Version
% \author{Sa\'ul A. Blanco and Esfandiar Haghverdi}
% \email{D}
% \date{January 3, 2025}
\makeatletter
\renewcommand\@date{{%
  \vspace{-2cm}%
  \large\centering
  \begin{tabular}{@{}c@{}}
    Sa\'ul A. Blanco \\
    \normalsize sblancor@iu.edu
  \end{tabular}%
  \quad and\quad
  \begin{tabular}{@{}c@{}}
    Esfandiar Haghverdi \\
    \normalsize ehaghver@iu.edu
  \end{tabular}

%  \bigskip

 Luddy School of Infomatics, Computing, and Engineering\\
 Indiana University, Bloomington

%   \bigskip

\today

\vspace{-1cm}
%  \date{February 12, 2025}
}}
\makeatother

\maketitle

\abstract{\vspace{-0.5cm} Motivated by the definition of Freiman homomorphism, we explore the possibilities of formulating some basic notions and techniques of additive combinatorics in a categorical language. We show that additive sets and Freiman homomorphisms form a category and we study several limit and colimit constructions in this, and in an interesting subcategory of this category. Moreover, we study the additive structure of these (co)limit objects using additive doubling constant. We relate this category to that of finite sets and mappings, and that of abelian groups and group homomorphisms. We show that the  Konyagin \& Lev result on the existence of universal ambient groups is an instance of adjunction.}

\section{Introduction}

The connections between category theory and combinatorial structures have been studied in different contexts. For example, Leinester~\cite{L08} defined and studied analogs of the Euler characteristic for categories. Moreover, the classic notion of M\"obius function and M\"obius inversion have been carried over to category theory~\cite{CLL80,L12, L75, S04}. There is also a well-established relationship between combinatorial objects, such as posets, and Hopf algebras (see, for example,~\cite{ABS06,AM10, E96, JR78, MR11}). Some work has also been done looking at inequalities that hold between certain homsets in some categories ~\cite{I91,KW03}. 

In this paper we give a categorical view of mathematical objects and some constructions associated with additive combinatorics. In essence, additive combinatorics studies properties of additive sets, and in particular, how large they can get with respect to  addition operation of their ambient group. We will briefly explain the basic notions in the next section. 

More precisely, we construct two categories, that we call $\FR$ and $\FR^0$. The former category uses additive sets as objects and Freiman homomorphisms as morphisms, and the latter utilizes additive sets that contain the additive identities of their ambient  groups as objects and Freiman homomorphisms that preserve this identity as morphisms. An important construction in additive combinatorics is that of the universal ambient group. We show that this construction can be seen as an instance of left adjoint construction. 

We now establish the basic additive combinatorics background that we need for the rest of this paper

\section{Additive Combinatorics}
We will start by a brief introduction to the main ideas and some of the results in additive combinatorics that we will be using in the rest of this paper. The main reference for this part of the paper is the classic book \cite{TV06} by Tao and Vu, see also~\cite{Zhao23} for a more recent source. 

Additive combinatorics is a branch of combinatorial number theory that studies the additive
structure of non-empty finite subsets of abelian groups. For example, consider the subset
$A = \{1,2,3\}$ of the abelian group $\mZ$ of integers. Then, $A + A = \{2,3,4,5,6\}$ and hence
$|A+A| = 5$ which is bigger than $|A|=3$. On the other hand, if we take the subset $B = \{0,1,2\}$ of the abelian group $\mZ_3$ of integers modulo $3$, we notice that $|B+B| = |B|$. 

Let us use a common measure of additive structure, namely the {\em doubling constant} $\sigma$ defined as 
$\sigma[A] = \frac{|A+A|}{|A|}$. Then, for the examples above $\sigma[A] = 5/3$ and $\sigma[B] = 1$. The closer the  doubling constant $\sigma[A]$ of a set $A$ is to 1, the more additive structure the set $A$ has. So in our case, $B$ has more additive structure than $A$, in fact $B$ is a group itself so it has maximal additive structure. A natural question arises as to whether the doubling constant could be arbitrarily large and the answer is affirmative, let $A =\{1,2,2^2,\ldots,2^{N-1}\}$, that is, a geometric progression, then one can verify that $\sigma[A] = (N+1)/2$. This shows that there are sets that can have minimal additive structure, that is, large $\sigma[A]$. 

On the other hand, we can generalize our first example above, that is , $A = \{1,2,3\}$ to an arithmetic progression  
$A = \{a,a+r,a+2r,\ldots, a+(N-1)r\}$ and observe that $\sigma[A] = 2 - 1/N$, so doubling constant can get as close to 2 as we wish. Additive combinatorics 
aims to characterize the sets by the range of their doubling constant. A fundamental result in the field is \emph{Freiman's theorem}, which we now state. By a \emph{generalized arithmetic progression}, one means a set of the form
\[
\{a_0+n_1a_1+\cdots+ n_da_d:0\leq n_i< N_i, \text{ with }1\leq i\leq d\}, 
\] where $a_0,\ldots,a_d,n_1,\ldots, n_d\in\mathbb{Z}$. Freiman's theorem states that if $A\subset\mathbb{Z}$ has a bounded doubling constant, namely $\sigma[A]\leq K$, then $A$ is contained in a generalized arithmetic progression  of dimension at most 
$d(K)$ and size at most $f(K)|A|$, where $d(K)$ and 
$f(K)$ are constants depending only on $K$, (see \cite[Chapter 7]{Zhao23}). 

The notion of Freiman homomorphism is a fundamental concept in additive combinatorics, comparable to the notion of group homomorphism for groups but more relaxed. Freiman homomorphisms play the same role for additive sets as group homomorphisms do for abelian groups (see \cite[Section 5.3]{TV06}). In particular, Freiman isomorphic additive sets have identical additive structure as measured by doubling constant, difference constant, and additive energy (see \cite{TV06} for the definition of these measures). 

Freiman homomorphism plays a critical role in the proof of Freiman’s theorem and the embedding of additive sets into universal ambient groups~\cite{TV06}.
We were motivated by such observations to organize additive sets and Freiman homomorphisms into a category and study this category for its own sake.
We have tried to make the exposition self-contained for both category theorists and mathematicians working in additive combinatorics and thus will assume just the basic definitions in each field.

\section{Category Theory in Combinatorics}
In this section, we motivate categorical thinking and its formulation in areas related to combinatorics.

In graph theory, one can construct the category of graphs, with graphs as objects and graph homomorphisms as morphisms. Graph constructions such as the \emph{tensor product} can be seen as the categorical product in the category of graphs and graph homomorphisms. In this context, statements such as Hedetniemi’s conjecture about the chromatic number of the product of two graphs $ G, H $, $ \chi(G \times H) \leq \min\{\chi(G), \chi(H)\} $, have a categorical interpretation in terms of morphisms. Moreover, pullbacks in the category of graphs have been studied extensively, as shown in~\cite{Brown08, Hell79, HellN04, Waller76}. Similarly, for hypergraphs, categorical methods have been explored in~\cite{DW80}.

Posets, including lattices, can be viewed as examples of \emph{thin} or \emph{posetal categories}, where each homset contains at most a unique morphism. Moreover, since a common theme in combinatorics is to ``add up  over all possible `decompositions' of an object," the language of bialgebras, where one has a product and a co-product, is very natural. Furthermore, Hopf algebras have found a special place in combinatorics, due in part to Aguiar, Bergeron, and Sotile's~\cite{ABS06} result establishing that the algebra of quasi-symmetric functions is a terminal object in the category of the so-called \emph{combinatorial Hopf algebras}. In their monograph, Aguiar and Mahajan~\cite{AM10} explore different connections between algebraic, combinatorial, and categorical objects.

Another significant development is Joyal’s theory of combinatorial species~\cite{J81} (also see~\cite{Bergeron97}). A \emph{combinatorial species} is defined as an endo-functor on the category of finite sets and bijections. The theory of combinatorial species is used for counting labeled structures through exponential generating functions. Furthermore, the theory can also be used for the enumeration of unlabeled versions of these labeled structures. Indeed, for any bijection from a set $ A $ onto itself, one obtains a corresponding bijection from the set of structures with label set $ A $ to itself. Thus, one has an action of the symmetric group on $ A $ that acts on these structures. The orbits of this action correspond to the unlabeled structures. In his seminal paper~\cite{J81}, Joyal includes a proof from ``THE BOOK" \cite{AZ08} of Cayley's formula that the number of labeled trees on $n$ vertices is $n^{n-2}$ utilizing combinatorial species.

 Whitehead~\cite{Whitehead1, Whitehead2, Whitehead3} established that for classical homotopy theory, simplicial sets are an adequate alternative to the category of topological spaces. More formally, there is a \emph{Quillen equivalence} between the category of topological spaces and the category of simplicial sets. Further connections between CW-complexes and their corresponding face posets were established, bridging algebraic topology and combinatorics (see, for example, Wachs lecture notes~\cite{W06}).

In addition, several combinatorial concepts have been extended to the categorical framework. For instance, the Euler characteristic and M\"obius function have been adapted to categories~\cite{L08, L12}. Similarly, Ramsey's theorem has been generalized to certain categories, as explored in~\cite{Graham72}. This extension of Ramsey’s theorem was instrumental in proving a conjecture of Rota stating that if $k,r,l$ are non-negative integers, and $\mathbb{F}_q$  a field with $q$ elements, then there is a positive integer $N=N(q,k,r,l)$ (in other words, $N$ depends only on $q,k,r,l$) that satisfies the following property: If $V$ is a vector space over $\mathbb{F}_q$ of dimension at least $N$, and if all $k$-dimensional subspaces of $V$ are divided into $r$ classes, then there exists some $l$-dimensional subspace with all of its $k$-dimensional subspaces in a single class. Rota's conjecture can be seen as a vector-space analogue of the classic Ramsey's theorem from the 1930s.

\section{Freiman Category}

We will continue with some basic definitions from additive combinatorics, the main reference here is the excellent book by Tao and Vu~\cite{TV06}. For category theory definitions, the reader can consult~\cite{BW90,Mac98}.

An {\em additive set} $A$ is a finite non-empty subset of an abelian ambient group $G$. We shall use $(A,G)$ to denote an additive set. Note that $A$ is simply a set and is not required to be a subgroup or a coset of a subgroup. We simply need to be able to add the elements of the set $A$.

\begin{defn}[Freiman Homomorphism]
{\em 
Let $k\geq 1$ and let $(A,Z)$ and $(B,W)$ be additive sets. A {\em Freiman homomorphism
$\phi : (A,Z) \to (B,W)$ of order $k$} is a map $\phi : A \to B$ such that
$a_1+a_2+\cdots + a_k = a'_1+a'_2+\cdots + a'_k$ implies $\phi(a_1)+\phi(a_2) + \cdots + \phi(a_k) = \phi(a'_1)+\phi(a'_2) + \cdots + \phi(a'_k)$, for all $a_1, a_2, \ldots, a_k, a'_1, a'_2, \ldots , a'_k \in A$. If furthermore $\phi$ has an inverse which is a Freiman $k$-homomorphism from $(B,W)$ to $(A,Z)$, then $\phi$ is a {\em Freiman $k$-isomorphism}. In that case, we say that $(A,Z)$ and $(B,W)$ are Freiman isomorphic of order $k$.
}
\end{defn}

Note that a Freiman homomorphism of order $k$ will also be a Freiman homomorphism of all orders $k'<k$. Any map from $A$ to $B$ is a Freiman homomorphism of order 1, which is why many authors require $k\geq 2$ in the definition of such a homomorphism. Also, any group homomorphism is clearly a Freiman homomorphism of any order $k \geq 1$.

Here are some more examples. The sets $A = \{0,1,2\}$ and $B=\{10,13,16\}$ as subsets of $\mZ$ are Freiman isomorphic for all $k$ via $\phi:A\to B$ given by $\phi(0)=10,\phi(1)=13,\phi(2)=16$. On the other hand, $C = \{0,1,4\}$ and $D = \{0,1,3\}$ are Freiman 3-isomorphic, for example, via $\phi:C\to D$ given by $\phi(0)=0,\phi(1)=1,\phi(4)=3$ but $C$ and $D$ are not 4-isomorphic. In fact, the only Freiman 4-homomorphisms from $C$ to $D$ are obtained by mapping every element in $C$ to the same element in $D$. 

It is easy to see that the identity map on any additive set is a Freiman homomorphism of any order $k$ and given Freiman homomorphisms $f: (A,Z) \to (B,W)$ and $g: (B,W) \to (C,V)$ of order $k$,
their function composition $g\circ f : (A,Z) \to (C,V)$ is a Freiman homomorphism of order $k$ as well. 
We can, therefore organize additive sets and Freiman homomorphisms of a fixed order $k$ between them  into a category defined below.
 
\begin{defn}[Freiman Category]
{\em Given an integer $k \geq 1$ we define the {\em Freiman category of order $k$}, denoted by
$\FR$, as the category whose objects are additive sets, that is a pair of the form $(A,Z)$ where $A$ is a finite non-empty subset of the additive group $Z$. The morphisms in $\FR$ are Freiman homomorphisms of order $k$.
}
\end{defn}

To explore categorical limits and colimits, we will also be working in a subcategory of $\FR$, which we define formally below. An additive set need not contain 0 (the additive identity of the ambient group), and thus a Freiman homomorphism is not required to preserve it. However, more structure is needed in order to consider categorical constructions. Indeed, we consider a subcategory of $\FR$ whose objects have more structure than additive sets do. We properly define this subcategory below.

\begin{defn}
{\em We say that an additive set $(A,Z)$ is \textit{normalized} if $0_Z\in A$, where $0_Z$ denotes the additive identity of $Z$. If $(B,W)$ is also normalized, and $f:(A,Z)\to(B,W)$ is a Freiman homomorphism, we say that $f$ {\em preserves} the group identity if  $f(0_Z)=0_W$, where $0_W$ is the additive identity of $W$. By abuse of notation, we shall drop the subscripts that indicate the underlying group, since this will be clear from the context. We denote by $\FR^0$ the subcategory of $\FR$ where the objects are normalized additive sets and the morphisms are Freiman homomorphisms of order $k$ that preserve the group identity. We call $\FR^0$ the \textit{Normalized Freiman Category} of order $k$. Note that $\FR^0$ is not a full subcategory of $\FR$, for example the translation map does not preserve the group identity.
}
\end{defn}

\textbf{Initial and Terminal Objects}--It is easy to see that $(\{0\},Z)$ (for any group $Z$), is both an initial and a  terminal object in $\FR^0$. In $\FR$, any object of the form $(\{*\},Z)$ (* is any element of $Z$), is a {\it weak} initial object (that is to say, given any other additive set $(B,W)$, there is, a not necessarily unique, map from $(\{*\},Z)$ to $(B,W)$), but it is a terminal object.

$\FR$ does not have an initial object as the only possible candidate for such an object
would be $(\emptyset,Z)$ (as in the category of sets and functions). This fails as  additive sets are defined to be non-empty.

We now prove some basic properties of morphisms in $\FR$.

\begin{prop}[Iso, Mono, Epi]
Let $f: (A,Z) \to (B,W)$ be a morphism in $\FR$ or $\FR^0$, then

\begin{enumerate}
\item[(a)] $f$ is an isomorphism of order $k$ iff it is a bijection  and $f^{-1}$ is also a Freiman homomorphism of order $k$. %\footnote{Note that this agrees with the definition of a Freiman isomorphism in \cite{TV06}.}
\item[(b)] $f$ is a monomorphism of order $k$ iff it is injective.
\item[(c)] $f$ is an epimorphism of order $k$ iff it is surjective.
\end{enumerate}
\end{prop}

\begin{proof}{}
The proof of (a) follows easily from the definition of an isomorphism in a category.
The proofs of parts (b) and (c) are quite straightforward following the arguments for the category of sets and functions. Indeed, for (b) suppose that $f$ is injective and consider two morphisms $g_1,g_2:(C,V)\to (A,Z)$ such that $f\circ g_1=f\circ g_2$. Then for any $c\in C$, it follows that $f(g_1(c))=f(g_2(c))$, and since $f$ is injective $g_1(c)=g_2(c)$. Thus, $g_1=g_2$ and $f$ is a monomorphism. Conversely, suppose that $f$ is a monomorphism and let $a,a'\in A$ satisfy $f(a)=f(a')$. Consider $g_1,g_2:(\{c\},V)\to(A,Z)$ given by $g_1(c)=a$ and $g_2(c)=a'$. Then $f\circ g_1=f\circ g_2$, and since $f$ is a monomorphism, $g_1=g_2$. Therefore $a=a'$ proving that $f$ is injective. For $\FR^0$ one uses the set $\{0,c\}$ and the proof goes through similarly. The proof of (c) is similar.   
\end{proof}

Note that a bijective Freiman homomorphism $\phi$ might not necessarily be a Freiman isomorphism since the inverse map $\phi^{-1}$ may not be a Freiman homomorphism. For example, consider the sets $A=\{0,1,4\}$ and $B=\{0,2,4\}$ along with $\phi: A\to B$ given by $\phi(0)=0,\phi(1)=2$, and $\phi(4)=4$. Note that $\phi$ is a Freiman 2-homomorphism from $A$ to $B$ but its inverse is not a Freiman 2-homomorphism, for example $2+2 = 4+0$ but $\phi^{-1}(2)+\phi^{-1}(2) = 1+1 \neq 4+0 = \phi^{-1}(4) + \phi^{-1}(0)$.

%What about EM factorization?

As mentioned earlier in this section, the key feature of a Freiman isomorphism is the preservation of the additive structure. 
Thus, for instance, a Freiman homomorphism of order 2 will not increase the size of $|A|$ or $|A + A|$ (see the theorem below for a precise statement). 

\begin{thm}[{\cite[Lemma 5.23]{TV06}}] 
Let $(A,G)$ be an additive set and $\phi : (A,G) \to (\phi(A),H)$ be a surjective Freiman homomorphism of order $k$. Then,

$$|\epsilon_1\phi(A_1) + \cdots + \epsilon_k \phi(A_k)| \leq |\epsilon_1A_1 + \cdots + \epsilon_kA_k|$$
where $A_1,\ldots,A_k$ are non-empty subsets of $A$ and $\epsilon_i = \pm1$. If $\phi$ is a Freiman isomorphim of order $k$, then we may replace inequality with equality. In particular, if $A$ and $B$ are Freiman isomorphic of order $k$, then
$$|lB-mB| = |lA-mA| \, whenever \, l,m\geq 0,\, l+m \leq k.$$
\end{thm}

% we shall also explore how epi and mono behave with respect to additive structure.

We now describe some categorical constructions in the category $\FR$. 

\begin{prop}[Product]\label{pr}
Given two objects $(A,Z)$ and $(B,W)$ in $\FR$ their categorical product is given by 
$(A\times B, Z \oplus W)$. Here $\times$ is the Cartesian product of sets and 
$\oplus$ is the direct sum (categorical biproduct) of abelian groups. 
\end{prop}

\begin{proof}{}
The projection map $pr_1 : (A\times B, Z \oplus W) \to (A, Z)$ is given by the restriction of projection map from $Z \oplus W$ to $Z$ (a group homomorphism) to the set $A\times B$. Similarly, for $pr_2$. It can be  easily verified that given any Freiman homomorphisms  $f : (C,V) \to (A,Z)$ and $g : (C,V) \to (B,W)$, there exists a unique mediating morphism $h: (C,V) \to  (A\times B, Z \oplus W)$ such that $pr_1 \circ h = f$ and $pr_2 \circ h = g$. Indeed, $h(c) = (f(c),g(c)).$

We need to show that the mediating morphism $h$ is a Freiman homomorphism.
Indeed, given $c_1,\dots,c_k$ and $c'_1,\dots, c'_k$ such that
$c_1+c_2+\dots + c_k = c'_1+c'_2+\dots + c'_k$, 
\begin{eqnarray*}
    h(c_1) + \dots + h(c_k) & = & (f(c_1),g(c_1)) + \dots + (f(c_k),g(c_k))\\ 
 & = & (f(c_1) + \dots + f(c_k), g(c_1) + \dots + g(c_k))\\
 & = & (f(c'_1) + \dots + f(c'_k), g(c'_1) + \dots + g(c'_k))\\
 & & \text{$f$ and $g$ are Freiman homomorphisms}\\
 & = & h(c'_1) + \dots + h(c'_k),
\end{eqnarray*} as desired.
\end{proof}

The argument in the proof above can be easily extended to show that $\FR$ has finite products.
The nullary product is the terminal object.

Similarly, one can show that $\FR^0$ has finite products defined in exactly the same way. In this case, 
$(0,0)$ is the group identity which is in $A\times B$. It can also be easily checked that all homomorphisms preserve the zero elements.

\begin{prop}[Coproduct]\label{copr}
Given two objects $(A,Z)$ and $(B,W)$ in $\FR^0$, let $\widehat{A} = A\times \{0\}$ and $\widehat{B} = 
\{0\}\times B$.
The categorical coproduct of $(A,Z)$ and $(B,W)$ is given  
by $(\widehat{A} \cup \widehat{B}, Z\oplus W)$. 
\end{prop}

\begin{proof}{}
The injection map $in_1 : (A,Z) \to (\widehat{A} \cup \widehat{B}, Z\oplus W)$ sends $a \in A$ to $(a,0)$. Similarly, $in_2 : (B,W) \to (\widehat{A} \cup \widehat{B}, Z\oplus W)$ sends $b \in B$ to $(0,b)$. Both $in_1$ and $in_2$ are Freiman $k$-homomorphisms as they are restrictions of group homomorphisms to the sets $A$ and $B$, respectively.

Given any other object $(C,V)$ in $\FR^0$ and Freiman homomorphisms $f: (A,Z) \to (C,V)$, and $g: (B,W) \to (C,V)$, there exists a unique mediating morphism $h : (\widehat{A} \cup \widehat{B}, Z\oplus W) \to (C,V)$ defined by 
$$h(e) = \begin{cases}
    f(a) & \text{if $e = (a,0)$,}\\
    g(b) & \text{if $e = (0,b)$.}\\
    %0 & \text{if $e = (0,0)$}
\end{cases}$$ 

such that $h \circ in_1 = f$ and $h \circ in_2 = g$. Note that $h$ is well-defined since $f(0)=0=g(0)$. 

We need to show that the mediating morphism $h$ is a Freiman homomorphism.

Indeed, Let $e_1, \dots, e_k$ and $e'_1,\dots, e'_k$ be in $\widehat{A} \cup \widehat{B}$
such that $$e_1 + \dots + e_k = e'_1 +  \dots + e'_k.$$ This identity holds if and only if 

$$a_1 + \dots + a_{k_1} = a'_1 + \dots + a'_{k_1}$$ 
and 
$$b_1 + \dots + b_{k_2} =  b'_1 + \dots + b'_{k_2}$$

for some $a_1, \dots, a_{k_1}, a'_1, \dots, a'_{k_1} \in A$, 
$b_1, \dots, b_{k_2}, b'_1,\dots, b'_{k_2} \in B$, and some $k_1,k_2$ such that $k_1+k_2 = k$.

Then
\begin{eqnarray*}
    h(e_1) + \dots + h(e_k) & = & f(a_1)+ \dots + f(a_{k_1}) + g(b_1) + \dots + g(b_{k_2})\\ 
 & = & f(a'_1) + \dots + f(a'_{k_1}) +  g(b'_1) + \dots + g(b'_{k_2})\\
 & & \text{$f$ and $g$ are Freiman homomorphisms}\\
 & = & h(e'_1) + \dots + h(e'_k),
\end{eqnarray*} as needed. 
\end{proof}

\begin{rem}
Unfortunately, the same construction does not yield a coproduct in the category $\FR$. The argument in the proof above can be easily extended to show that $\FR^0$ has finite coproducts where the nullary coproduct is the initial object. 
\end{rem}

We now turn our attention to pullbacks and pushouts, both of which exist in $\FR^0$.

\begin{prop}[Pullback]\label{pback}
The category $\FR^0$ has pullbacks. 
\end{prop}

\begin{proof}{}
Given normalized additive sets $(A,Z),(B,W)$, and $(C,V)$ and Freiman homomorphisms $f : (A,Z) \to (C,V)$, and $g: (B,W) \to (C,V)$ in $\FR^0$, the \emph{fiber product} is defined as $A\times_C  B = \{(a,b)\,|\, f(a) = g(b)\}$. Then the pullback of $f$ and $g$ is defined by $((A\times_C B, Z\oplus W),p_1,p_2)$ where $p_1(a,b)=a$, and $p_2(a,b) = b$. Note that $f\circ p_1 = g \circ p_2$, and $(0,0) \in A\times_C B$. Given normalized $(D,U)$ and Freiman homomorphisms $d_1:(D,U)\to (A,Z)$ and $d_2:(D,U)\to(B,W)$ with $f\circ d_1 = g\circ d_2$, there exists a unique $h$, namely $h(x)=(d_1(x),d_2(x))$, such that $d_i=p_i\circ h$, for $i=1,2$. Note that, by definition $f\circ d_1 (x) = g\circ d_2(x)$ for all $x \in D$ and thus $h(x)$ is in 
$A \times_C B$. 

Next we show that the mediating morphism $h$ is a Freiman homomorphism that preserves additive identities.
Indeed, let $a_1,\dots, a_k$ and $a'_1,\dots, a'_k$ be in $D$ such that
$a_1 + \dots + a_k = a'_1 + \dots + a'_k$. Then

\begin{eqnarray*}
    h(a_1) + \dots + h(a_k) & = & (d_1(a_1),d_2(a_1)) + \dots + (d_1(a_k),d_2(a_k))\\ 
 & = & (d_1(a_1) + \dots + d_1(a_k), d_2(a_1) + \dots + d_2(a_k))\\
 & = & (d_1(a'_1) + \dots + d_1(a'_k), d_2(a'_1) + \dots + d_2(a'_k))\\
 & & \text{$d_1$ and $d_2$ are Freiman homomorphisms}\\
 & = & (d_1(a'_1),d_2(a'_1)) + \dots + (d_1(a'_k),d_2(a'_k))\\
 & = & h(a'_1) + \dots + h(a'_k)
\end{eqnarray*}

\begin{displaymath}
    \xymatrix{(D,U) \ar@{.>}[drr]_h \ar@/_/[ddrr]_{d_2} \ar@/^/[drrrr]^{d_1}   &&&&\\
&  &  (A\times_C B, Z\oplus W) \ar[d]^{p_2} \ar[rr]_{p_1} &
                 & (A,Z)  \ar[d]_f       \\
&  & (B,W) \ar[rr]^g &  & (C,V)                }
\end{displaymath}
Furthermore, since $h(x)=(d_1(x),d_2(x))$, $h$ preserves the additive identity. 
\end{proof}

Note that this construction will not yield pullbacks in $\FR$ as the fiber product of two sets may be empty.
Note, also that this is the only impediment to the existence of pullbaks in $\FR$.

\begin{prop}[Pushout]\label{pout}
The category $\FR^0$ has pushouts.
\end{prop}

\begin{proof}{} 
Given additive sets $(A,Z),(B,W)$, and $(C,V)$ and Freiman homomorphisms 
$f:(C,V)\to(A,Z)$ and $g:(C,V)\to(B,W)$,
we define the pushout of $f,g$ by 
$(((\widehat{A}\cup \widehat{B})/\sim,(Z\oplus W)/N),i_1,i_2)$. Here, $\widehat{A} = A\times \{0\}$ and $\widehat{B} = \{0\}\times B$, $N$ is the subgroup of 
$Z\oplus W$ generated by the set $\{(f(c),-g(c)),\,\, c\in C\}$, and $\sim$ is the smallest equivalence relation on $\widehat{A}\cup \widehat{B}$ such that $(a,b) \sim (a',b')$ if
$(a,b) - (a',b') = (f(c),-g(c))$, for some $c \in C$. 

The morphism $i_1: (A,Z) \to ((\widehat{A}\cup \widehat{B})/\sim, (Z\oplus W)/N)$ is defined by $i_1(a) = [(a,0)]$, and $i_2: (B,W) \to ((\widehat{A}\cup \widehat{B})/\sim, (Z\oplus W)/N)$ is defined by $i_2(b) = [(0,b)]$. Here, $[.]$ denotes the equivalence class under the relation $\sim$.

Given $(D,U)$ and Freiman homomorphisms $d_1: (A,Z)\to (D,U)$ and $d_2:(B,W)\to (D,U)$ with $d_1\circ f = d_2\circ g$, there exists a unique $h : ((\widehat{A}\cup \widehat{B})/\sim, (Z\oplus W)/N) \to (D,U)$, defined by 
$$h(e) = \begin{cases}
    d_1(a) & \text{if $e = [(a,0)]$,}\\
    d_2(b) & \text{if $e = [(0,b)]$.}\\
    %0 & \text{if $e = [(0,0)]$}
\end{cases}$$ 

such that $d_1 = h\circ i_1$ and
$d_2 = h\circ i_2$. Note that $h$ is well-define since $d_1(0)=d_2(0)=0$.

Next, we show that the mediating morphism $h$ is a Freiman homomorphism. Below we have omitted the brackets for the sake of readability. Indeed, let $e_1, \dots, e_k$ and $e'_1,\dots, e'_k$ be in $\widehat{A} \cup \widehat{B}$
such that $$e_1 + \dots + e_k = e'_1 +  \dots + e'_k.$$ This identity holds if and only if 

$$a_1 + \dots + a_{k_1} = a'_1 + \dots + a'_{k_1}$$ 
and 
$$b_1 + \dots + b_{k_2} =  b'_1 + \dots + b'_{k_2}$$ 

for some $a_1, \dots, a_{k_1}, a'_1, \dots, a'_{k_1} \in A$, 
$b_1, \dots, b_{k_2}, b'_1,\dots, b'_{k_2} \in B$, and some $k_1,k_2$ such that $k_1+k_2 = k$.

Then
\begin{eqnarray*}
    h(e_1) + \dots + h(e_k) & = & d_1(a_1)+ \dots + d_1(a_{k_1}) + d_2(b_1) + \dots + d_2(b_{k_2})\\ 
 & = & d_1(a'_1) + \dots + d_1(a'_{k_1}) +  d_2(b'_1) + \dots + d_2(b'_{k_2})\\
 & & \text{$d_1$ and $d_2$ are Freiman homomorphisms}\\
 & = & h(e'_1) + \dots + h(e'_k),
\end{eqnarray*} as needed. 

 \begin{displaymath}
     \xymatrix{(D,U)     &&&&\\
 &  &  ((\widehat{A}\cup \widehat{B})/\sim,(Z\oplus W)/N) \ar@{.>}[ull]_h   &
                  & (A,Z) \ar@/_/[ullll]_{d_1}  \ar[ll]_{i_1}       \\
 &  & (B,W)\ar[u]_{i_2}  \ar@/^/[uull]^{d_2}  &  & (C,V)   \ar[u]_f     \ar[ll]^g        }
\end{displaymath} 

\end{proof}

Unfortunately,  the same construction does not yield a pushout in the category $\FR$. 

\begin{prop}[Equalizers and coequalizers]\label{eq}
The category $\FR^0$ has equalizers and coequalizers.
\end{prop}

\begin{proof}{}
Given parallel morphisms $f,g : (A,Z) \to (B,W)$ their equalizer is given by
$((E,Z),e)$ where $E = \{a\in A \,|\, f(a) = g(a)\}$ and $e : E \to A$ is simply the inclusion map. Given another additive set $(C,V)$ and Freiman homomorphism $k:(C,V)\to(A,Z)$ such that $f\circ k=g\circ k$, then 
there exists a unique map $h:(C,V) \to (E,Z)$ defined as $h(c)=k(c)$ for all $c\in C$ such that $e\circ h = k$. Clearly, $h$ is a Freiman homomorphism.

\begin{displaymath}
    \xymatrix{  (E,Z) \ar[r]^e &  (A,Z)\ar@<1ex>[r]^f \ar[r]_g &  (B,W)  &\\
     (C,V) \ar[ru]^k \ar[u]^h&  &  & }
\end{displaymath}

The coequalizer of $f$ and $g$ is given by 
$((B/\sim,W/\im(f-g)),c)$ where $\sim$ is the smallest equivalence relation on $B$ generated by 
$\{(f(a),g(a))\mid a\in A\}\subseteq B\times B$ and $c$ is defined by $c(b)=[b]$, the equivalence class of $b$ with respect to $\sim$, for all $b \in B$. Given an additive set $(D,Y)$ and a Freiman homomorphism $k:(B,W)\to (D,Y)$ satisfying $k\circ f=k\circ g$, there exists a unique map $h:(B/\sim,W/\im(f-g))\to(D,Y)$ given by $h([b])=k(b)$ for all $[b] \in B/\sim$ such that
$h \circ c = k$. Since $k$ is a Freiman homomorphism, so is $h$. 

\begin{displaymath}
    \xymatrix{(A,Z)\ar@<1ex>[r]^f \ar[r]_g & (B,W) \ar[dr]_k \ar[r]^c   & \,\,\,\,(B/\sim,W/\im(f-g))\ar[d]^h & \\
      &  &  (D,Y)   & }
\end{displaymath}
\end{proof}

Note that $\FR$ does not have equalizers as $E = \{a\in A \,|\, f(a) = g(a)\}$ may be empty.
Similarly, the coequalizer construction above will not yield a coequalizer in $\FR$ as $\im(f-g)$ may not form a subgroup of $W$ (as $0$ may not be in $\im(f-g)$).

\begin{cor} Combining Propositions \ref{pr} and \ref{pback} (\ref{copr} and \ref{pout}), we have established that $\FR^0$ is a complete (co-complete) category.
\end{cor}

\begin{rem}
$\FR^0$ is not a \emph{Cartesian closed category} because it has a zero object and any category with a zero object and products is Cartesian closed if and only if it is equivalent to a category with only one object and the identity morphism.
\end{rem}

Normalized additive sets and zero preserving Freiman homomorphisms, $\FR^0$, interact very nicely with the limit and colimit objects, which shows that $\FR^0$ is more well-behaved than $\FR$ in this respect.

\subsection{Comparison to other categories}
Now that we have explored the categorical properties of $\FR^0$ and $\FR$, we should compare them to two closely related categories, namely the category \textbf{FSet} of finite sets and mappings, and the category \textbf{Ab} of abelian groups and group homomorphisms. 

There is a well-defined forgetful functor $U : \FR \to \textbf{FSet}$ sending an additive set $(A,Z)$ to the underlying set $A$ that forgets the ambient group $Z$. A Freiman homomorphism $f : (A,Z) \to (B,W)$ will be mapped to $U(f) : A \to B$, $U(f)=f$. It is easy to check that this relationship is functorial: $U(1_A) =1_A$ and for any $f:(A,Z) \to (B,W)$ and $g:(B,W) \to (C,V)$, $U(gf) = gf = U(g)U(f).$ 

Interestingly enough the correspondence between $\FR$ and \textbf{Ab} is not that straightforward. That is because not every Freiman homomorphism can be uniquely extended to a group homomorphism, which is needed to make sure that the functorial relationship works properly at the level of morphisms. However, there exists an object in additive category that is relevant in this context, the \textit{universal ambient group}, which we define below using the already established categorical notation above.

\begin{defn}
{\em 
Let $(A,Z)$ be an object in $\FR$. $Z$ is said to be a \emph{universal ambient group} for $A$ if for any other additive set $(B,W)$ and $f : (A,Z) \to (B,W)$, $f$ can be uniquely extended to a group homomorphism from $Z$ to $W$.
}
\end{defn}

We shall therefore limit ourselves to a subcategory of $\FR$ that we call $\UFR$ where the ``U'' stands for universal. The objects of category $\UFR$ are additive sets $(A,Z)$ where $Z$ is universal for $A$. The morphisms are the same as in $\FR$. Thus $\UFR$ is a full subcategory of $\FR$.

We can now proceed with the definition of another forgetful functor $V : \UFR \to \textbf{Ab}$ given by $V (A, Z) = Z$.  Thus $V$ forgets the underlying set $A$, and given $f : (A, Z) \to (B,W)$, $V(f) : Z \to W$ is the unique group homomorphism that extends $f$ and is guaranteed to exist by definition.

%So we can limit ourselves to finite abelian groups and group homomorphisms and then it works. Is there an adjointness here?

\begin{defn}
{\em
An additive group $Z'$ is a {\em universal ambient group} for $(A,Z)$ if there exists an object $(A',Z')$ which is $k$-isomorphic to $(A,Z)$ and such that 
$Z'$ is a universal ambient group for $A'$.
}
\end{defn}

\begin{thm}[Konyagin \& Lev, \cite{KL00}]
Fix $k\geq 2$, and let $(A,Z)$ be an additive set. Then there exists a universal group $Z'$ for $A$. Furthermore, if $A'$ is an embedding of $A$ inside this ambient group $Z'$, then $Z'$ is generated as a group by $A'$. In particular, $Z'$ is finitely generated.
\end{thm}

Consider the functor $F : \FR \to \FR$, where
$F(A) = (A',Z')$ with $Z' = A^*/\langle X\rangle$, $A^*$ is the free abelian group generated by some basis $\{e_a \,|\, a\in A\}$, and $\langle X\rangle$ is the subgroup of $A^*$ generated by the set $$X = \{e_{a_1}+\dots + e_{a_k} - e_{a'_1}- \dots - e_{a'_k}\,|\, a_1+\dots + a_k = a'_1 +\dots + a'_k, a_1,\dots,a_k,a'_1,\dots,a'_k \in A\}.$$ 

$A'$ is the image of $A$ under the maps $A \stt{e} A^* \stt{\pi} A^*/\langle X\rangle$, where $\pi e(a) = \pi(e_a)$. Here the first  map is the standard embedding map sending $a$ to the basis element $e_a$, and $\pi$ is the canonical quotient map. The additive set $(A',Z')$ as constructed is the universal ambient group for $(A,Z)$ (see the proof of~\cite[Theorem 5.39]{TV06}).

As for the morphisms, given $f: (A,Z) \to (B,W)$, $F(f) : (A',Z') \to (B',W')$ is defined by $F(f)(e_a + \langle X\rangle) =
e_{f(a)} + \langle f(X) \rangle$. Here $\langle f(X)\rangle$ is the subgroup of $B^*$ generated by  
\begin{eqnarray*}
f(X) & = & \{e_{f(a_1)} + \cdots + e_{f(a_k)} - e_{f(a'_1)} - \cdots - e_{f(a'_k)} \,|\, \\
& & f(a_1)+\cdots + f(a_k) = f(a'_1) +\cdots + f(a'_k), \\ & & a_1,\dots,a_k,a'_1,\dots,a'_k \in A\}.
\end{eqnarray*}

We have the following result.
\begin{prop}
The universal ambient group construction of Konyagin and Lev is an instance of adjunction. That is, the functor $F$ above is left adjoint to the identity functor $Id$. In other words, $F$ is naturally isomorphic to the identity functor via  a unique isomorphism.
\end{prop}

\begin{proof}{}
It is easy to see that for any Freiman $k$-homomorphism $f : (A,Z) \to (B,W)$, $F(f)$ is a Freiman $k$-homomorphism from $(A',Z') \to (B',W')$.

$F(1_{(A,Z)}
): (A',Z') \to (A',Z')$ is given by $F(1_{(A,Z)})(e_a +\langle X\rangle) = e_a + \langle X\rangle = 1_{F(A,Z)}$. Also, given $f: (A,Z) \to (B,W) $ and $g: (B,W) \to (C,V)$, $F(g f) (e_a + \langle X \rangle) = e_{gf(a)} + \langle gf(X)\rangle 
 = F(g) (e_{f(a)} + \langle f(X) \rangle) = F(g) (F(f) (e_a + \langle X \rangle))$ and hence $F$ is a functor. As for the adjointness, we show the following bijection:
 $$\FR(F(A,Z), (B,W)) \cong \FR((A,Z), Id((B,W)).$$
Given a Freiman $k$-homomorphism $f: (A',Z') \to (B,W)$ in the left hand set above, we define
$\theta(f) : (A,Z) \to (B,W)$ by $\theta(f)(a) = f(e_a + \langle X \rangle)$ which is clearly a Freiman $k$-homomorphism. On the other hand, given a Freiman $k$-homomorphism $g: (A,Z) \to (B,W) $ in the right hand set above, 
$\eta(g) : (A',Z') \to (B,W)$ is defined by $\eta(g)(e_a + \langle X \rangle) =  g(a)$. Note that $\eta(g)$ is also a Freiman $k$-homomorphism. Finally, it is an easy exercise to show that $\theta $ and $\eta$ above are inverses of each other.
\end{proof}

Inspired by the construction above we can define another adjunction. This will require a modification of the definition of an additive set. We define a new category $\FR^\infty$ where the objects are pairs $(A,Z)$ with $Z$ an abelian group and $A$ a non-empty subset of $Z$. Note that in contrast to an additive set, the set $A$ here does not have to be finite. Morphisms are as before, Freiman $k$-homomorphisms. Even though this takes us away from the realm of additive combinatorics, such infinite additive sets are the main object of study in the field of additive number theory. Let $\textbf{Ab}$ be the category of abelian groups and group homomorphisms. Define the functor $F : \FR^\infty \to \textbf{Ab}$ on objects by $F((A,Z)) = A^*/\langle X \rangle$ where all is as defined above. Given a morphism $(A,Z) \to (B,W)$, $F(f) : A^*/\langle X \rangle \to B^*/ \langle Y \rangle$ is defined as a group homomorphism determined by its action on the basis elements of $A^*/\langle X \rangle$ where $F(f)(e_a + \langle X \rangle) = e_{f(a)} + \langle f(X) \rangle$.  
It is easy to show that $F$ is a functor. Now define a functor $G : \textbf{Ab} \to \FR^\infty$ by $F(Z) = (Z,Z)$ and given a group homomorphism $f: Z \to W$, $F(f) =f$

\begin{prop}
Let $F$ and $G$ be two functors as defined above. Then $F$ is left adjoint to $G$.
\end{prop}

\begin{proof}{}  We show the following bijection holds:
 $$\textbf{Ab}(F(A,Z),W) \cong \FR^\infty((A,Z), G(W)).$$
Given a group homomorphism $f: A^*/\langle X \rangle \to W$ in the left hand set above, we define
$\theta(f) : (A,Z) \to (W,W)$ by $\theta(f)(a) = f(e_a + \langle X \rangle)$ which is clearly a Freiman $k$-homomorphism. On the other hand, given a Freiman $k$-homomorphism $g: (A,Z) \to (W,W) $ in the right hand set above, 
$\eta(g) : A^*/\langle X \rangle \to W$ is defined as the group homomorphism determined by $\eta(g)(e_a+ \langle X \rangle) = g(a)$. Finally, one can easily verify that $\theta$ and $\eta$ above are inverses of each other.
\end{proof}

Finally, we describe an alternative categorical description for the universal ambient group based on exercises 5.5.1 and 5.5.4 in \cite{TV06}. We will show that this construction also yields an instance of adjunction. 

Let $\mT = \mR/\mZ$ be the circle group. Given an additive set $(A,Z)$, consider the set $G = \FR((A,Z), (\mT,\mT))$ of Freiman $k$-homomorphisms from $(A,Z)$ to $(\mT,\mT)$.
This set is an abelian group where the sum is defined by addition. It can be viewed as a compact subgroup of a torus and as such we can take its \emph{Pontryagin dual}  $Z' = \hat{G}$ which consists of all continuous group homomorphisms from the abelian group $G$ to $\mT$. We define $A' \subseteq Z'$ as 
$A' = \{\hat{a} \,|\,  a \in A\}$ where $\hat{a} \in Z'$ is defined by $\hat{a}(\chi) = \chi(a)$ for all $\chi \in G$. One can show that $Z'$ is a universal ambient group for $A$. Or in other words, that  $(A',Z')$ is Freiman $k$-isomorphic to $(A,Z)$ and $Z'$ is a universal ambient group for $A'$. Now define the functor $F: \FR \to \FR$ by $F((A,Z)) = (A',Z')$ as defined above. Given $f : (A,Z) \to (B,W)$, $F(f): (A',Z') \to (B',W')$ is defined by $F(f)(\hat{a}) = \widehat{f(a)}$.

\begin{prop}
The universal ambient group construction above is an instance of adjunction. That is, the functor $F$ above is naturally isomorphic to the identity functor via  a unique isomorphism. 

\end{prop}
\begin{proof}{}
First we need to show that given a Freiman $k$-homomorphism $f: (A,Z) \to (B,W)$, $Ff : (A',Z') \to (B',W')$ as constructed above is a Freiman $k$-homomorphism. We let $G = \FR((A,Z),(\mT,\mT))$ and $G' = \FR((B,W),(\mT,\mT))$. 

Suppose $a_1 + \cdots + a_k = a'_1 + \cdots + a'_k$ for $a_1,\ldots, a_k, a'_1, \ldots a'_k \in A$, then $f(a_1) + \cdots + f(a_k) = f(a'_1) + \cdots + f(a'_k)$ as $f$ is a Freiman $k$-homomorphism. Moreover, $\chi'(f(a_1)) + \cdots + \chi'(f(a_k)) = \chi'(f(a'_1)) + \cdots + \chi'(f(a'_k))$ as $\chi' \in G'$ is a Freiman $k$-homomorphism.

Now, suppose $\hat{a_1} + \cdots +\hat{a_k} = \hat{a'_1} + \cdots + \hat{a'_k}$ for $\hat{a_1}, \ldots, \hat{a_k}, \hat{a'_1}\ldots, \hat{a'_k} \in A'$. For any $\chi' \in G'$,

\begin{eqnarray*}
    (Ff(\hat{a}_1) + \cdots + Ff(\hat{a}_k))(\chi') & = & \widehat{f(a_1)}(\chi') + \ldots + \widehat{f(a_k)}(\chi')\\
    & = & \chi'(f(a_1)) + \cdots + \chi'(f(a_k))\\
    & = & \chi'(f(a'_1)) + \cdots + \chi'(f(a'_k)), \, \text{using the identity above}\\
    & = & \widehat{f(a'_1)}(\chi') + \cdots + \widehat{f(a'_k)}(\chi')\\
    & = & (Ff(\hat{a'_1}) + \cdots + Ff(\hat{a'_k}))(\chi')
\end{eqnarray*}

Thus establishing the fact that $Ff$ is a Freiman $k$-homomorphism.

$F(1_{(A,Z)}): (A',Z') \to (A',Z')$ is given by $F(1_{(A,Z)})(\hat{a}) = \hat{a} = 1_{F(A,Z)}$. Also, given $f: (A,Z) \to (B,W) $ and $g: (B,W) \to (C,V)$, 
$F(gf)(\hat{a}) = \widehat{g(f(a))} = Fg(\widehat{f(a)}) = Fg(Ff (\widehat{a}))$. Hence $F$ is a functor. To prove adjointness, we show the following bijection:
 $$\FR(F(A,Z), (B,W)) \cong \FR((A,Z), Id((B,W)).$$
Given a Freiman $k$-homomorphism $f: (A',Z') \to (B,W)$ in the left hand set above, we define
$\theta(f) : (A,Z) \to (B,W)$ by $\theta(f)(a) = f(\hat{a})$ which is clearly a Freiman $k$-homomorphism. On the other hand, given a Freiman $k$-homomorphism $g: (A,Z) \to (B,W) $ in the right hand set above, 
$\eta(g) : (A',Z') \to (B,W)$ is defined by $\eta(g)(\hat{a}) = g(a)$. Note that $\eta (g)$ is also a Freiman $k$-homomorphism. 

Finally, it is an easy exercise to show that $\theta $ and $\eta$ above are inverses of each other.
\end{proof}

Motivated by the previous construction, we define a new adjunction for the category $\FR^\infty$ as follows. Consider the category $\textbf{LCA}$ of locally compact abelian groups where the objects are locally compact abelian groups and morphisms are group homomorphisms. Now consider the functor $F: (\FR^\infty)^{op} \to \textbf{LCA}$ defined by 
$$F((A,Z)) = \FR^\infty((A,Z), (\mT,\mT)).$$
Given objects $(A,Z)$ and $(B,W)$ and a map 
$f:(A,Z) \to (B,W)$ in $\FR^\infty$,
$$F(f) : \FR^\infty((B,W), (\mT,\mT)) \to \FR^\infty((A,Z), (\mT,\mT))$$
is given by $F(f)(g) = g\circ f$. 

Now define the functor $H: \textbf{LCA} \to (\FR^\infty)^{op}$ by 
$H(Z) = (\widehat{Z},\widehat{Z})$ where $\widehat{Z}$ is the Pontryagin dual of $Z$. Given $f : Z \to W$, $H(f) : (\widehat{W},\widehat{W}) \to (\widehat{Z},\widehat{Z})$ is defined by $H(f)(\chi) = \chi \circ f$ for any $\chi \in \widehat{W}$.

\begin{prop}
Let $F$ and $H$ be two functors as defined above. Then $H$ is left adjoint to $F$.
\end{prop}

\begin{proof}{}  We show the following bijection holds:
 $$\FR^\infty((B,W),H(Z)) \cong \textbf{LCA}(Z, F(B,W)).$$
Given a Freiman $k$-homomorphism $f: (B,W) \to (\widehat{Z},\widehat{Z})$ in the left hand set above, we define
$\theta(f) : Z \to \FR^\infty ((B,W),(\mT,\mT))$ by $\theta(f)(z)(b) = f(b)(z)$ for $z\in Z$ and $b \in \FR^\infty ((B,W),(\mT,\mT))$. $\theta(f)$ is a continuous group homomorphism.

Moreover, given a continuous group homomorphism $g: Z \to \FR^\infty ((B,W),(\mT,\mT))$ in the right hand set above, 
$\eta(g) : (B,W) \to (\widehat{Z},\widehat{Z})$ is defined as $\eta(g)(b)(z) = g(z)(b)$ for $z\in Z$ and $b\in B$. Note that $\eta(g)$ is a Freiman $k$-homomorphism. Finally, it is an easy exercise to show that $\theta$ and $\eta$ above are inverses of each other.
\end{proof}

\section{The Additive Structure of (co)Limits}{\label{s:addstructure}}

We now explore the additive structure of the sets that underlie the additive sets in various limit and colimit constructions in $\FR^0$. Even though disjoint union of additive sets does not occur in any limit or colimit construction we have added it to our list as well. We first recall the definition of the doubling constant mentioned in the introduction. 

\begin{defn}
For an additive set $(A,Z)$, $\sigma(A):=\frac{|A+A|}{|A|}$ is called the \em{doubling constant} of $(A,Z)$
\end{defn}

\begin{prop}[Additive Structure]
Let $(A,Z)$ and $(B,W)$ be additive sets. Then

\begin{enumerate} 
\item[(i)] $\sigma[\{0\}] =1$.
\item[(ii)] $ \sigma[A\times B] = \sigma[A]\sigma[B].$
\item[(iii)] $1\leq \sigma[A \times_C B] \leq \sigma[A\times B]$, where $(A,Z)$ and $(B,W)$ are normalized additive sets, and $A \times_C B = \{(a,b)\,|\, f(a) = g(b)\}$ for $f: A \to C$ and $g: B \to C$ that preserve additive identities.
\item[(iv)] $1 \leq \sigma[A/\sim]\leq \sigma[A]$, where $\sim$ is an equivalence relation on $A$.
\item[(v)]  $1 \leq \sigma[E] \leq \sigma[A]$, where $E = \{a \in A \,|\, f(a) = g(a)\}$ for $(A,Z)$ and $(B,W)$ objects in $\FR^0$ and $f,g: (A,Z) \to (B,W)$ morphisms in $\FR^0$.

\item[(vi)] $\sigma(A\uplus B) \leq \sigma(A)+\sigma(B)+\frac{|A+B|}{|A| + |B|}$, where $A$ and $B$ have the same ambient group, that is in this case only $Z=W$.

\end{enumerate}

\end{prop}

\begin{proof}{}
\begin{itemize}
    \item[(i)] Follows from $|\{0\}+\{0\}|=1=|\{0\}|$.

\item[(ii)] Notice that $(A\times B)+(A\times B)=\{(a+a',b+b')\mid a,a'\in A\text{ and }b,b'\in B\}$. Hence,
\begin{align*}
   \frac{|(A\times B)+(A\times B)|}{|A\times B|}&=\frac{|\{(a+a',b+b')\mid a,a'\in A\text{ and }b,b'\in B\}|}{|A\times B|}\\
   &=\frac{(|A+A|)(|B+B|)}{|A\times B|}\\
   &=\frac{|A+A|}{|A|}\frac{|B+B|}{|B|}\\
   &=\sigma(A)\sigma(B),
\end{align*} as desired. 

\item[(iii)] Follows from the fact that $1 \leq |A\times_C B| \leq |A\times B|$. The lower bound happens when $\{(a,b) \in A\times B \,|\, f(a) = g(b)\} = \{(0,0)\}$ and the upper bound happens when $\{(a,b) \in A\times B \,|\, f(a) = g(b)\} = A\times B$.

\item[(iv)] The set of equivalence classes could be as small as a single class and as large as the size of $A$.

\item[(v)] The set $E$ could be as small as $\{0\}$ and as large as $A$.

\item[(vi)] Note  that

\begin{align*}
    \sigma(A\uplus B)&=\frac{|A\uplus B+A\uplus B|}{|A\uplus B|}\\
    &\leq\frac{|A+A|+|A+B|+|B+B|}{|A\uplus B|}\\
    &= \frac{|A+A|}{|A\uplus B|}+\frac{|A+B|}{|A\uplus B|}+\frac{|B+B|}{|A\uplus B|}\\
    &\leq \frac{|A+A|}{|A|}+\frac{|A+B|}{|A\uplus B|}+\frac{|B+B|}{|B|}\\
    &=\sigma(A)+\sigma(B)+\frac{|A+B|}{|A\uplus B|}\\
    &=\sigma(A)+\sigma(B)+\frac{|A+B|}{|A|+|B|},
\end{align*} as desired
\end{itemize}

\end{proof}

\section{Concluding remarks}

In this paper, we study two categories formed by additive sets and Freiman homomorphisms, $\FR$ and $\FR^0$. We provide several limit and colimit constructions including pullbacks, pushouts, equalizers, and coequalizers. We furthermore establish that Konyagin and Lev's universal ambient group construction and another construction due to Tao and Vu using Pontryagin duality are instances of adjunctions.

The work in this paper is taking the first baby steps towards a more in-depth categorical study of additive combinatorics. Such an study, we believe, could facilitate a technology transfer from  fields like algebraic topology and algebraic geometry were categorical thinking has proven so effective and successful in overcoming major obstacles into additive combinatorics. 
The technology transfer from these field into additive combinatorics could open up the use of (co)homological techniques in proving results in additive and arithmetic combinatorics, as well as leading to new discoveries.

In the other direction, we believe that many ideas and techniques in additive combinatorics could lead to interesting results in the combinatorial studies of locally finite (all homsets are finite sets) categories.

As a more immediate next step, we are currently focusing on formulating and proving Freiman's structure theorem in a categorical framework. It seems that this would require enriching the world of objects and morphisms with suitably designed quantitative measures so that we can keep track of the complexity (size) of various structures that occur in this theorem and various others, for example, the Freiman-Ruzsa polynomial conjecture.

\section*{Acknowledgments} The authors thank Omer Cantor for several helpful comments about an earlier draft of the paper.

\end{document}